\documentclass[10pt,twoside]{article}
\usepackage{graphicx}
\usepackage{amsmath}
\usepackage{Latex-document}
\usepackage{latexsym}

\title{\bf Analysis of Some Singular Solutions  \vskip -1mm in Fluid Dynamics\thanks{Supported in part by
grants from the Research Grants Council of Hong Kong Special
Administrative Region CUHK 4219/99D and CUHK 4279/00P.}\vskip 6mm}

\newcommand{\nc}{\newcommand}
\nc{\ptl}{\partial}
\newtheorem{theorem}{Theorem}[section]

\input{amssym.def}
\nc{\R}{\Bbb{R}} \nc{\RR}{{\R}^2} \nc{\RRR}{{\R}^2_+}
\nc{\Rr}{{\R}^3} \nc{\M}{\mathcal{M}}

\author{Zhouping Xin\vspace*{-0.5cm}\thanks{The Institute of Mathematical Sciences \& Department of Mathematics,
The Chinese University of Hong Kong, Shatin, N.T., Hong Kong.
E-mail: zpxin@ims.cuhk.edu.hk}}
\date{\vspace{-8mm}}

\begin{document}

\maketitle

\thispagestyle{first} \setcounter{page}{851}

\begin{abstract}\vskip 3mm
Studies on singular flows in which either the velocity fields or
the vorticity fields change dramatically on small regions are of
considerable interests in both the mathematical theory and
applications.  Important examples of such flows include supersonic
shock waves, boundary layers, and motions of vortex sheets, whose
studies pose many outstanding challenges in both theoretical and
numerical analysis.  The aim of this talk is to discuss some of
the key issues in studying such flows and to present some recent
progress.  First we deal with a supersonic flow past a perturbed
cone, and prove the global existence of a shock wave for the
stationary supersonic gas flow past an infinite curved and
symmetric cone.  For a general perturbed cone, a local existence
theory for both steady and unsteady is also established. We then
present a result on global existence and uniqueness of weak
solutions to the 2-D Prandtl's system for unsteady boundary
layers. Finally, we will discuss some new results on the analysis
of the vortex sheets motions which include the existence of 2-D
vortex sheets with reflection symmetry; and no energy
concentration for steady 3-D axisymmetric vortex sheets.

\vskip 4.5mm

\noindent {\bf 2000 Mathematics Subject Classification:} 35L70,
35L65, 76N15.

\noindent {\bf Keywords and Phrases:} Singular flows, Shock waves, Vortox sheets, Prandtl's system, Boundary
layers.
\end{abstract}

\vskip 12mm

\section{Introduction} \label{section 1}\setzero
\vskip-5mm \hspace{5mm }

Many physically interesting phenomena involve evolutions of
singular flows whose velocity fields or vorticity fields change
dramatically.  Shock waves, vortex sheets, and boundary layer are
some of the well-known examples of such singular flows which
provide better approximations for significant parts of the flow
fields near the physical boundaries, in the mixing layers and
trailing wakes, etc., at high Reynolds numbers. The better
understanding of the dynamics of such singular flows is the key in
the analysis of general fluid flows governed by the well-known
Euler (or Navier-Stokes) systems or their variants for both
compressible and incompressible fluids, and has been one of the
main focuses for applied analysists for decades.  Substantial
progress has been made in the past in studying singular flows by
either rigorous analysis, or numerical simulations, or asymptotic
methods [18].  In particular, a rather complete theory exists for
the 1-D shock wave problems, and both the theoretical
understanding and numerical methods for 2-D smooth incompressible
flows are quite satisfactory.  Yet there are still many important
issues to be settled, such as, the global existence of weak
solutions and asymptotic structures of various approximate
solutions generated by either physical considerations or numerical
methods for such singular flows.

One of the fundamental problems in the mathematical theory of
shock waves for hyperbolic conservation laws is the well-posedness
of the multi-dimensional gas-dynamical shock waves, for which the
celebrated Glimm's method does not apply.  Most of the previous
studies along this line deal with either short time structural
stability of basic fronts or asymptotic analysis and numerical
simulations of dynamics of shock fronts.  Due to the great
complexity and the lack of understanding, it is reasonable for one
to begin with some of the physically relevant wave patterns where
a lot of experimental data, numerical simulations, and asymptotic
results are available.  One of the basic model for such studies is
the supersonic flow past a pointed body [4], which is one of the
fundamental problems in gas dynamics.  We first study the
stationary supersonic gas flow past an infinite curved and
symmetric cone.  The flow is governed by the potential equation as
well as the boundary conditions on the shock and the surface of
the body.  This problem has been studied extensively by either
physical experiments or numerical simulations.  The rigorous
analysis starts with the work of Courant and Friedrichs in [4],
where they show that if a supersonic flow hits a circular cone
with axis being parallel to the velocity of the upstream flow and
the vertex angle being less than a critical value, then there
appears a circular conical shock attached at the tip of the cone,
and the flow field between the shock front and the surface of the
body can be determined by solving a boundary value problem of a
system of ordinary differential equations.  The local existence of
supersonic flow past a pointed body has been established recently
[1].  Our interest is on the structure of the global solution to
this problem.  We show that the solution to this problem exists
globally in the whole space with a pointed curved shock attached
at the tip of the cone and tends to a self-similar solution [2].
Our analysis is based on a global uniform weighted energy estimate
for the linearized problem.  The method we developed in [2] seems
to be quite effective for other multi-dimensional problems.
Indeed, similar approach can be used to study unsteady supersonic
flows past a curved body for both potential flows and the full
Euler system, and we obtain the local existence of shock waves in
these cases [3].

Another challenging problem in the mathematical theory of
fluid-dynamics is the theoretical foundation of the Prandtl's
boundary layer theory [14].  In the presence of physical
boundaries, the solutions to the invisicid Euler system cannot be
the uniform asymptotic ansatz of the corresponding Navier-Stokes
system for large Reynolds number due to the discrepancies between
the no-slip boundary conditions for the Navier-Stokes system and
the slip boundary condition for the Euler system. Indeed, the
physical boundaries will creat vorticity and there is a thin layer
(called boundary layer) in which the leading order approximation
of the flow velocity is governed by the Prandtl's system for the
boundary layers (see (3.1) in Section 3). There are extensive
literatures on the theoretical, numerical, and experimental
aspects of the Prandtl's boundary layer theory [14]. Yet very
little rigorous theory exists for the dynamical boundary layer
behavior of the Navier-Stokes solutions for both compressible and
incompressible fluids.  One of the main difficulties is the
well-posedness theory in some standard H\"{o}lder or Sobolov
spaces for the initial-boundary value problems for the Prandtl's
system for boundary layers which is a severely degenerate
parabolic-elliptic system, for which the only known existence
results are proved locally in the analytic class in [15], except
the series of important works of Oleinik who dealt with a class of
monotonic data [12].  Indeed, Oleinik considered a plane unsteady
flow of viscous incompressible fluid in the presence of an
arbitrary injection and removal of the fluid across the boundary.
Under the monotonicity assumption (see (3.5) in Section 3),
Oleinik proved the well-posedness of local classical solutions to
the Prandtl's system [12].  One of the open questions posed in
[12] is to prove the global well-posedness of solution for the
Prandtl's system under suitable conditions. Recently, we establish
a global existence and uniqueness of weak solutions to the 2-D
Prandtl's system for unsteady boundary layers in the class
considered by Oleinik provided that the pressure is favorable.
This is achieved by introducing a viscous splitting method and new
weighted total variation estimate [16,17].  See Section 3 and
[16,17] for more details.

Finally, we turn to the motion of vortex sheets, which corresponds
to a singular inviscid flow where the vorticity field is zero
except on lower-dimensional surfaces, the sheets, and can be
characterized as inviscid flows with finite local energy and with
vorticity fields being finite Radon measures.  The study on the
existence and structure of solutions for the inviscid Euler system
for incompressible fluids with data in such class is of
fundamental importance both physically and mathematically.
Physically, vortex sheets can be used to model important flows
such as high Reynolds number shear layers, and have many
engineering applications. Mathematically, the evolution of a
vortex sheets gives a classical example of ill-posed problem in
the sense of Hadamard, a curvature singularity developes in finite
time, and the nature of the solution past singularity formation is
of great interest to know. This gives rise to many interesting yet
challenging problems. Some of these are:  Is there a
well-posedness theory of classical weak solutions to the inviscid
Euler system with general vortex sheets initial data? What are the
structures of the approximate solutions to vortex sheets motions
generated by either Navier-Stokes solutions or pratical numerical
methods (such as particle method)?  Can vorticity concentration
and energy defects occur dynamically? etc.. Despite the importance
of these problems and past intensive effort in rigorous
mathematical analysis, these problems are far from being solved.
Better understanding has been achieved before singularity
formation in the analytical setting, and studies on global weak
solutions and their approximations start with the important works
of Diperna-Majda [6].  Delort observes that no vorticity
concentration implies that a weak limit in $L^2$ of an approximate
solution sequences is in fact a classical weak solution to the
two-dimensional Euler system, and thus proved the first existence
of global (in time) classical weak solution to the 2-D
incompressible Euler equations with vortex sheets initial data
provided that the initial vorticity is of distinguished sign [5].
Similar ideas have been used to study the convergence of
approximate solutions generated by either viscous regularization
[11] or partical methods [8,9] for vortex sheets with one sign
vorticity.  In the case that vorticity may change sign, the vortex
sheets motion becomes extremely complex after singularity
formation. Indeed, many important features of irregular flows seem
to be connected with interactions and intertwining of regions of
both positive and negative vorticities. Here we consider a
mirror-symmetric flow which allows interactions but excludes
intertwining of regions of distinguish vorticity, prove that there
is no vorticity concentrations for approximate solution of such
flow, and thus give the first global (in time) existence of vortex
sheets motion with two-sign vorticity.  We also consider the 3-D
axisymmetric vortex motions.  It is well-known that for smooth
flows, the analysis for the axisymmetric 3-D Euler system without
swirls is almost identitial to that of 2-D incompressible Euler
system. However, this parallelness breaks down for the vortex
sheets motions.  Indeed, we will show that there exist no energy
concentrations in approximate solutions to vortex sheets motion
with one-sign vorticity for steady axisymmetric 3-D Euler system
without swirls [7].  Some partial results on unsteady axisymmetric
vortex sheets motion will also be discussed.

\section{The supersonic flow past a pointed body} \label{section 2}
\setzero\vskip-5mm \hspace{5mm}

A projectile moving in the air with supersonic speed, is governed
by the inviscid compressible Euler systems
\begin{eqnarray}
\left\{
\begin{array}{ll}
  \ptl_t\,\rho + div (\rho u) = 0, & \\
  \ptl_t (\rho u) + div (\rho u \otimes u) + \nabla p = 0, & \\
\end{array}
\right.
\end{eqnarray}
\noindent where $\rho, u = (u_1, u_2, u_3)$ and $p$ stand for the
density, the velocity and the pressure respectively.  We will only
treat the polytropic gases so that $p = p(\rho) = A \rho^\gamma$
with gas constant $A > 0$ and $1 < \gamma < 3$, $\gamma$ being the
adiabatic exponent.

Suppose that there is a uniform supersonic flow $(u_1, u_2, u_3) =
(0, 0, q_0)$ with constant density $\rho_0 > 0$ which comes from
negative infinity.  Then the flow can be described by the steady
Euler system.  If we assume further that the flow is irrotational,
so that one can introduce a potential function $\Phi$ such that $u
= \nabla \Phi$.  Then the Bernoulli's law implies that $\rho =
h^{-1} \left( \frac{1}{2} q^2_0 + h (\rho_0) - \frac{1}{2}
{|\nabla \Phi|}^2 \right) \equiv H(\nabla \Phi)$, where $h(\rho)$
is the specific enthalpy defined by $h^{'} (\rho) =
\frac{p^{'}(\rho)}{\rho}$.  In this case, (2.1) is reduced to a
second order quasilinear equation
\begin{equation}
div (H (\nabla \Phi) \nabla \Phi) = 0,
\end{equation}
\noindent which can be verified to be strictly hyperbolic with
respect to $x_3$ if $\ptl_3 \Phi > c$ with $c$ being the sound
speed given by $c^2 (\rho) = p^{'} (\rho)$.  The flow hits a point
body, whose surface is denoted by $m (x_1, x_2, x_3) = 0$.  Since
no flow can cross the boundary, the natural boundary condition is
\begin{equation}
u \cdot \nabla m \equiv \nabla \Phi \cdot \nabla m = 0 \quad {\rm
on} \quad m (x_1, x_2, x_3) = 0.
\end{equation}

If the vertex angle of the tangential cone of the pointed body is
less than a critical value, it is then expected that a shock front
is attached at the tip of the pointed body.  Denote by $\mu (x_1,
x_2, x_3) = 0$ the equation of the shock front, then the
Rankine-Hugoniot conditions become
\begin{equation}
\nabla \mu \cdot [H(\nabla \Phi) \nabla \Phi] = 0, \quad \Phi
\quad {\rm is\ continuous}, \quad {\rm on} \quad \mu (x_1, x_2,
x_3) = 0.
\end{equation}

Our aim is to find a solution to this free boundary value problem,
(2.2)--(2.4).  When the pointed body is small perturbation of a
circular cone, the local existence of solution to the problem
(2.2)--(2.4) has been established in [1].  Our main goal is to
establish a global solution.  However, such a global shock wave
might not exist in general for arbitrary pointed body due to the
possibility of development of new shock waves in the large. Thus
we assume further that the pointed body is a curved and symmetric
cone.  In this case, it will be more convenient to rewrite the
problem (2.2)--(2.4) in terms of polar coordinates $(r, \theta,
z)$ with $r = \sqrt{x^2_1 + x^2_2}$ and $z = x_3$. Assume that the
tip of the pointed body locates at the origin, the equation of the
surface of the body is $r = b(z)$ with $b(0)= 0$, and the equation
of the shock front is $r = S(z)$ with $S(0) = 0$.  Set $\Phi =
q_0\,z + \varphi(r,z)$.  Then (2.2)-- (2.4) become
\begin{equation}
\left( {\left(q_0 + \ptl_z \phi \right)}^2 - c^2 \right)
\ptl^2_{zz} \phi + \left( {\left( \ptl_r\,\phi \right)}^2 - c^2
\right) \ptl^2_{rr}\,\varphi + 2 \ptl_r\,\varphi (q_0 +
\ptl_z\,\varphi) \ptl^2_{rz}\,\varphi -
\frac{c^2}{r}\ptl_r\,\varphi = 0,
\end{equation}
\begin{equation}
-(q_0 + \ptl_z\,\varphi) b^{'}(z) + \ptl_r\,\varphi = 0 \quad {\rm
on} \quad r=b(z),
\end{equation}
\begin{equation}
-[(q_0 + \ptl_z\,\varphi)H] S^{'}(z) + [\ptl_r\,\varphi\,H] = 0
\quad {\rm on} \quad r=S(z).
\end{equation}
\noindent Moreover, the potential $\varphi (r,z)$ is continuous on
the shock, so it should satisfy $\varphi(S(z),z) = 0$.  Then in
[2], we have shown that the problem (2.5)--(2.7) has a globally
defined solution as summarized in the following theorem:
\begin{theorem}
{Assume that a curved and symmetric cone is given such that $b(0)=
0, \ b^{'}(0)=b_0, \ b^{(k)} (0)=0, \ 2 \leq k \leq k_1$, and
\begin{equation}
|z^k \frac{d^k}{dz^k}(b(z) - b_0 z)| \leq \varepsilon_0 \quad {\rm
for} \quad 0 \leq k \leq k_2, \ z>0
\end{equation}
\noindent with $k_1$ and $k_2$ being some suitable integers.
Suppose that a supersonic polytropic flow parallel to the $z$-axis
comes from negative infinite with velocity $q_0$, and density
$\rho_0>0$. Then for suitably small $\varepsilon_0$, $b_0$ and
$q^{-1}_0$, the boundary value problem (2.5)--(2.7) admits a
global weak entropy solution with a pointed shock front attached
at the origin. Moreover, the location of the shock front and the
flow field between the shock and the surface of the body tend to
the corresponding ones for the flow past the unperturbed circular
cone $r=b_0 z$ with the rate $z^{-\frac{1}{4}}$.}
\end{theorem}

It should be noted that there are no other discontinuities in our
solution besides the main shock.  Since the deviation of the
surface of the body from that of a circular cone is sufficiently
small (see (2.8)), any possible compression of the flow will be
absorbed by the main shock.  This is the mechanism to prevent the
formation of any new shocks inside the flow field caused by the
perturbation of the body.  In particular, our results demonstrate
that self-similar solution with a strong shock is structurally
stable in a global sense.  Indeed, the key element in the proof of
Theorem 2.1 is to establish some global uniform weighted energy
estimates for the linearized problem of (2.5)--(2.7) around the
self-similar solution with a strong shock obtained when the
pointed body is a circular cone.  This is achieved by a deliberate
choice of multipliers which must satisfy a system of ordinary
differential inequalities with complicated coefficients due to the
structure of the background self-similar solution and the
requirement of obtaining global estimates independent of $z$ for
the potential function and its derivatives on the boundary as well
as its interior of a domain.

When the projectile changes its speed, or it confronts some airstream, then the flow around the projectile will be
time-dependent.  Thus, we also consider the unsteady supersonic flow past a pointed body.  Although our analysis
applies to more general case [3], here we will only present the result for two-dimensional polytropic, unsteady,
and irrotational flow past a curved wedge.  For simplicity in presentation, we also assume that both the wedge and
the perturbed incoming flow from infinity are symmetric about $x_1$-axis.  Let $x_2 = b(x_1)$ with $b(0) = 0$ be
the equation of the wedge, and $x_2 = S(t, x_1)$ with $S(t, x_1=0)=0$ be the equation of the shock front.  Then in
terms of the velocity potential function $\phi$ (so that $u_1=\ptl_1\,\phi$ and $u_2=\ptl_2\,\phi$), we are
looking for solutions to the following initial boundary value problem
$$ \ptl_t \left( H(-\phi_t + \frac{1}{2}
{|\nabla \phi|}^2 ) \right) + \Sigma^2_{i=1} \ptl_{x_i} \left( \ptl_{x_i}\,\phi\,H (-\ptl_t\,\phi - \frac{1}{2}
{|\nabla \phi|}^2 ) \right) = 0, $$
\begin{equation}
  \hspace*{5cm} t>0, x_2>b_1(x_1),
\end{equation}
\begin{equation}
\nabla \phi \cdot (b^{'} (x_1), -1) = 0 \quad {\rm on} \quad x_2 =
b(x_1),
\end{equation}
\begin{equation}
[H] \ptl_t\,S + [H\,\ptl_{x_1}\,\varphi] \ptl_{x_1}\,S -
[H\,\ptl_{x_2}\,\phi] = 0 \quad {\rm on} \quad x_2 = S(t, x_1),
\end{equation}
\begin{equation}
S(0, x_1) = S_0(x_1), \ \phi(0, x_1, x_2) = \phi_0 (x_1, x_2), \
\ptl_t\,\phi (0, x_1, x_2) = 0,
\end{equation}
\noindent where $S_0 (x_1)$ and $\phi_0 (x_1, x_2)$ are suitably
small perturbations of the corresponding shock location and
potential function respectively for the steady flow.  Then we have
the following local existence result [3].

\begin{theorem}
{There exist positive constant $\delta_1$ and $\delta_2$ and functions $S(t, x_1)$ and $\phi(t, x_1, x_2)$ defined
on the regions $\{(t,x_1)|0<t<\delta_1, 0<x_1<\delta_2\}$ and $\{(t, x_1, x_2)|0$ $<t<\delta_1, 0<x_1<\delta_2,
b(x_1)<x_2<S(t,x_1)\}$ respectively such that $(S(t,x_1)$, \linebreak $\varphi(t,x_1,x_2))$ solves the problem
(2.9)--(2.12).}
\end{theorem}

\section{Prandtl's system for boundary layers} \label{section 3}
\setzero\vskip-5mm \hspace{5mm}

Consider a plane unsteady flow of viscous incompressible fluid in
the presence of an arbitrary injection and removal of the fluid
across the boundaries.  In this case, the corresponding Prandtl's
system takes the form
\begin{equation}
\ptl_t\,u + u\,\ptl_x\,u + v\,\ptl_y\,u + \ptl_x\,p =
\nu\,\ptl^2_{yy}\,u, \quad \ptl_x\,u + \ptl_y\,v = 0,
\end{equation}
\noindent in the region, $R=\{(x,y,t)|0 \leq x \leq L, \ 0 \leq y
< + \infty, \ 0<t<T\}$, where $\nu, L$ and $T$ are positive
constants. The initial and boundary conditions can be imposed as
\begin{eqnarray}
\left\{
\begin{array}{ll}
u|_{t=0} = u_0(x,y), \ u|_{x=0} = u_1 (y,t), \ u|_{y=0} = 0, \
v|_{y=0} = v_0(x,t), \\
\lim_{y \rightarrow + \infty} u(x,y,t) = U(x,t),
\end{array}
\right.
\end{eqnarray}
\noindent with $U=U(x,t)$ given as is determined by the
corresponding Euler flow.  The pressure $p=p(x,t)$ in (3.1) is
determined by the Bernoulli's law:  $\ptl_t\,U + U\,\ptl_x\,U +
\ptl_x\,p=0$.

It follows from the physical ground that one may assume that
\begin{equation}
U(x,t)>0, \ u_0(x,y)>0, \ u_1(y,t)>0, \ {\rm and} \ v_0(x,t) \leq
0.
\end{equation}

Due to the degeneracy in the Prandtl's system (3.1), the problem
of well-posedness theory of solutions to the problem (3.1)--(3.3)
in the standard H\"{o}lder space or Sobolov space is quite
difficult.  In a series of important works by Oleinik and her
coauthors [12], they studied this problem under the additional
assumption that the data are monotonic in the sense that
\begin{equation}
\ptl_y\,u_0(x,y)>0, \ {\rm and} \ \ptl_y\,u_1(y,t)>0,
\end{equation}
\noindent and prove that there exists a unique local classical
smooth solution to the initial-boundary value problem (3.1)--(3.3)
provided that the data are monotonic in the sense of (3.4). Here
by local we mean that $T$ is small if $L$ is given and fixed, and
$T$ is arbitrary if $L$ is small.  One of the open problem in [12]
is:  What are the conditions ensuring the global in time existence
and uniqueness of solutions to the problem (3.1)--(3.4) for
arbitrarily given $L$?  In [16,17], we study such problem and
establish the global (in time) existence and uniqueness of weak
solutions to the initial-boundary value problem (3.1)--(3.4) in
the case that the pressure is favorable, i.e.,
\begin{equation}
\ptl_x\,p(x,t) \leq 0, \ t>0, \ 0<x<L.
\end{equation}

More precisely, we have ([16,17]):
\begin{theorem}
{Consider the initial-boundary value problem for the 2-D Prandtl's
system, (3.1)--(3.2).  Assume that the initial and boundary data
satisfy the constraints (3.3), (3.4) and (3.5). Then there exists
a unique global bounded weak solutions to the initial-boundary
value problem (3.1)--(3.2).  Furthermore, these solutions are
Lipschitz continuous in both space and time.}
\end{theorem}

We remark here that the condition (3.5), that the pressure is
favorable, is exactly what fluid-dynamists believe for the
stability of a laminar boundary layer.  This is also consistent
with the case of stationary flows [12].  In the case of pressure
adverse, i.e., $\frac{\ptl p}{\ptl x}>0$, separation of boundary
layer may occur, so one would not expect the long time existence
of solution to (3.1)--(3.3).  Finally, we remark that in the case
that (3.5) fails as in many pratical physical situations, short
time existence of regular solution is still expected, which has
not been established yet.

\section{Vortex sheets motions} \label{section 4}
\setzero\vskip-5mm \hspace{5mm}

Two-dimensional vortex sheets motion corresponds to an inviscid
flow whose vorticity is zero except on one-dimensional curves, the
sheets.  Thus, it is governed by the following Cauchy problem
\begin{eqnarray}
\left\{
\begin{array}{ll}
\ptl_t\,w + u \cdot \nabla w = 0, \ u=K \ast w, \\
w(t=0, x) = w_0(x)\,\epsilon\,\M(\RR) \cap H^{-1}_{loc} (\RR),
\end{array}
\right.
\end{eqnarray}
\noindent where $u=(u_1,u_2)$ is the velocity field, $w=
\nabla^\perp \cdot u$ is the vorticity, and $K$ is the Biot-Sawart
kernel, and $\M(\RR)$ denotes the space of finite Radon measures
in $\RR$.  Let $(u^\varepsilon, w^\varepsilon) (w^\varepsilon
\equiv \nabla^\perp \cdot u^\varepsilon)$ be a sequence of
approximate solution to (4.1) with the following bounds
\begin{equation}
\sup_{0 \leq t \leq T} \ \int_{\RR} |w^2 (x,t)| dx \leq C_1 (T), \
\sup_{0 \leq t \leq T} \ \int_{|x| \leq R} {|u^\varepsilon
(x,t)|}^2 dx \leq C_2 (T,R).
\end{equation}

Then there exist $u \epsilon L^\infty \left( 0,T,L^2_{loc} (\RR)
\right)$ and $w \epsilon L^\infty \left( 0,T,\M(\RR) \right)$ with
$w=\nabla^\perp \cdot u$ such that
\begin{equation}
w^\varepsilon \rightharpoonup w \ {\rm in} \ \M \left( [0,T]
\times \RR \right), \ u^\varepsilon \rightharpoonup u \ {\rm in} \
L^2_{loc} \left( [0,T] \times \RR \right)
\end{equation}

Two main questions arise:  Is $(u,w)$ a classical weak solution to
problem (4.1)?  Does either vorticity concentration (i.e.
${\overline{\lim}}_{\varepsilon, r \rightarrow 0^+} \ \int_{B_r}
|w^\varepsilon (x,t)| dx > 0$) or energy defects (
${\overline{\lim}}_{\varepsilon \rightarrow 0} \ \int_{|x| \leq R}
{|u^\varepsilon (x,t)|}^2 dx > \int_{|x| \leq R} {|u(x,t)|}^2 dx$
for some $t>0, \ R>0$) occur?  It is clear from the structure of
the 2-D Euler system that no energy defects implies strong
$L^2$-convergence of the velocity field and thus the existence of
the classical weak solution to (4.1).  A less obvious fact is that
no vorticity concentration also implies the weak limit $(u,w)$
being a classical solution to (4.1), which follows from the
vorticity formulation of (4.1) as observed by Delort [5].  Then
the convergence to a classical weak solution to the Cauchy problem
(4.1) is proved for approximate solutions generated by either
regularizing the initial data [5], or Navier-Stokes approximation
[11], or vortex blob methods [8], or point vortex methods [9]
provided that the initial vorticity is of distinguished sign.  To
study the corresponding issues for flows where interactions of
regions of both positive and negative vorticities are allowed, in
[10], we study a 2-D mirror-symmetric vortex sheets motion whose
vorticity is a integrable perturbation of a non-negative
mirror-symmetric radon measure.  Here a Radon measure $\mu$ is
said to be non-negative mirror-symmetric (NMS) if $\mu|_{\RRR}
\geq 0$ and $\mu$ is odd with respect to $x_1=0$.  Then the main
results in [10] can be summarized as
\begin{theorem}
{Assume that $w_0 \equiv \mu_1 + \mu_2$ such that
$\mu_1\,\epsilon\,\M_c(\RR) \cap\,H^{-1}_{loc} (\RR) \cap NMS$ and
$\mu_2\,\epsilon\,L^1_c (\RR)$.  Then there exists a global (in
time) classical solution $(u,w) \epsilon\,L^\infty \left( 0,T,L^2
(\RR) \right) \otimes L^\infty \left( 0,T,\M (\RR) \right)$ to the
initial value problem (4.1).  Furthermore, this weak solution can
be obtained as a limit of either a sequence of smooth inviscid
solutions or a sequence of solutions to the Navier-Stokes system.}
\end{theorem}

It should be noted that this is the only result of existence of
classical weak solution to (4.1) involving vorticities with
different signs.  This is proved by showing $\int^T_0 \
\sup_{x_0\,\epsilon\,\RR} \ \int_{B(x_0,\delta)} \ |w^\varepsilon
(y,t)| dy\,dt \rightarrow 0$ as $\delta \rightarrow 0^+$ uniformly
in $\varepsilon$, i.e., no vorticity concentration occur any
where.  Theorem 4.1 also indicates that interactions of regions of
positive and negative vorticities without interwining may not
cause concentration in vorticity.  There remain many important
open problems for the 2-D vortex sheets motion such as the
existence of classical weak solution to (4.1) for general vortex
sheets initial data, and whether energy defects occur dynamically
even in the case of one-sign vorticity.

Finally, we consider 3-D axisymmetric vortex sheets motions.  In
cylindrical coordinate, $(r,\theta,z)$, axisymmetric solutions of
3-D Euler system have the form
\begin{equation}
u (x,t) = u^r (r,z,t) e_r + u^\theta (r,z,t) \ e_\theta + u^z
(r,z,t) \ e_z, \ p (x,t) = p (r,z,t)
\end{equation}
\noindent where $e_r = (\cos \theta, \sin \theta, 0)$, $e_\theta=
(-\sin \theta, \cos \theta,0)$, and $e_z=(0,0,1)$.  The
axisymmetric flow is said without swirls if $u^\theta \equiv 0$.
In this case, the vorticity field is given by $w= \nabla \times
u=w^\theta\,e^\theta$ with $w^\theta=\ptl_r\,u^z - \ptl_z\,u^r$,
and ${\tilde{D}}_t (r^{-1}\,w^\theta)=0$ with ${\tilde{D}}_t =
\ptl_t + u^r\,\ptl_r +u^z\,\ptl_z$.  Thus, for smooth data, the
theory for 3-D axisymmetric Euler system without swirls is almost
parallel to that of 2-D Euler equation.  However, this similarity
breaks down for vortex sheets motions.  Indeed, in sharp contrast
to the 2-D case, we show in [7] that there are no energy defects
for a suitable sequence of approximate solutions for 3-D steady
axisymmetric vortex sheets motion with one-signed vorticity.
Precisely, we have
\begin{theorem}
{Let $(u^\varepsilon, p^\varepsilon)$ be smooth axisymmetric
solution to the 3-D steady Euler system: $(u^\varepsilon \cdot
\nabla) u^\varepsilon + \nabla p^\varepsilon=f^\varepsilon$,
$div\,u^\varepsilon=0$, $x\,\epsilon\,\Rr$, for some given
axisymmetric function $f^\varepsilon$ with $f^\varepsilon
\rightharpoonup f$ weakly in $L^1 (\Rr)$.  Suppose further that
\begin{equation}
{(w^\varepsilon)}^\theta \geq 0, \ \sup_\varepsilon \ \int_{\Rr} \
|w^\varepsilon| dx < + \infty, \ \sup_\varepsilon \ \int_{\Rr} \
{|u^\varepsilon|}^2 dx < + \infty
\end{equation}
\noindent where $w^\varepsilon=\nabla \times u^\varepsilon =
{(w^\varepsilon)}^\theta \,e_\theta$.  Let $u$ be the weak limit
of $u^\varepsilon$ in $L^2 (\Rr)$.  Then $u$ is a classical weak
solution to $(u \cdot \nabla) u + \nabla p=f$, $div\,u=0$,
$x\,\epsilon\,\Rr$.  Moreover, there exists a subsequence
$\{u^{\varepsilon_j} \}$ of $\{u^\varepsilon\}$ such that
$u^{\varepsilon_j}$ converges to $u$ strongly in $L^2 (\Rr)$.}
\end{theorem}

The proof of this theorem is based on a shielding  method and the
following fact which is valid for both steady and unsteady flows
[7].

\begin{theorem}
{Let $\{u^\varepsilon\}$ be a sequence of approximate solutions
for 3-D axisymmetric Euler system with general vortex-sheets data
generated by either smoothing the initial data or Navier-Stokes
approximations.  Let $A$ be the support of the defect measure
associated with $u^\varepsilon\,\epsilon\,L^2 ([0,T] \times \Rr)$.
If $A \neq \phi$, then $A\,\cap \{ (x,t)\epsilon\,\Rr \times
[0,T]| \ r>0 \} \neq \phi$.}
\end{theorem}

Theorem 4.3 implies that if there are no energy defects away from
the symmetry axis, then strong $L^2$-convergence takes place.  It
remains to study whether energy defects occur for unsteady
axisymmetric 3-D Euler equations.

\label{lastpage}

\end{document}